\documentclass[12pt,leqno]{article}

\newcounter{conjecture}
\newcounter{remark}
\newcounter{corollary}

\newcommand{\eqnsection}{
   \renewcommand{\theequation}{\thesection.\arabic{equation}}
   \makeatletter
   \csname @addtoreset\endcsname{equation}{section}
   \makeatother}
\newtheorem{theorem}{Theorem}
\newtheorem{lemma}{Lemma}

\def \be{\begin{equation}}
\def \ee{\end{equation}}
\def \bt{\begin{theorem}}
\def \et{\end{theorem}}
\def \bea{\begin{eqnarray}}
\def \eea{\end{eqnarray}}
\def \bas{\begin{eqnarray*}}
\def \eas{\end{eqnarray*}}

\def \al{\alpha}
\def \bb{\beta}
\def \ga{\gamma}
\def \Ga{\Gamma}
\def \de{\delta}

\def \ep{\epsilon}

\def \la{\lambda}
\def \La{\Lambda}

\def \om{\omega}
\def \Om{\Omega}

\def \si{\sigma}

\def \th{\theta}

\def \ze{\zeta}

\def \ff{\infty}
\def \wh{\widehat}
\def \wt{\widetilde}

\def \rar{\rightarrow}

\def \DD{{\cal D}}

\def \({\left(}
\def \){\right)}

\def \lc{\left\{}
\def \rc{\right\}}

\def \nn{\nonumber}

\def \bc{\begin{center} }
\def \ec{\end{center} }
\def \bs{\begin{slide} }
\def \es{\end{slide} }

\def\square{{\vcenter{\vbox{\hrule height.3pt
            \hbox{\vrule width.3pt height5pt \kern5pt
               \vrule width.3pt}
            \hrule height.3pt}}}}
\def\qed{{\hfill $\square$ \bigskip}}

\def\ol{\overline}

\eqnsection
\begin{document}

\title{An almost sure invariance principle for the range of  planar
random walks}

\author{Richard F. Bass\thanks{Research partially supported by
NSF grant
\#DMS-0244737}
\, \,and  \, Jay Rosen\thanks {Research partially supported by
grants from the NSF  and from PSC-CUNY.}}



\maketitle

\bibliographystyle{amsplain}

\begin{abstract} For a symmetric random walk in $Z^2$ with
$2+\de$ moments, we represent
$|\mathcal{R}(n)|$, the cardinality of  
the range, in terms of  an expansion involving the
renormalized intersection local times of a Brownian motion.  We
show that  for each $k\geq 1$
\[ (\log n)^k \left[ \frac{1}{n} |\mathcal{R}(n)| +\sum_{j=1}^k
(-1)^j (\textstyle{\frac1{2\pi}}\log n +c_X)^{-j}
\ga_{ j,n}\right]\to 0, \qquad a.s.
\] where $W_t$ is a Brownian motion, $W^{(n)}_t=W_{nt}/\sqrt
n$, 
$\ga_{ j,n}$ is the renormalized intersection local time at time 1
for
$W^{(n)}$, and
$c_X$ is a constant depending on the distribution of the random
walk.
\end{abstract}

\section{Introduction}

Let $S_{ n}=X_{ 1}+\cdots+X_{ n}$ be a random walk in $Z^2$,
where
$X_{ 1},X_{2},\ldots$ are symmetric i.i.d. vectors in $Z^2$. We
assume that the
$X_{i}$ have $2+\de$ moments for some $\de>0$ and covariance
matrix equal to the identity. We assume futher that the random
walk
$S_{ n}$ is strongly aperiodic in the sense of Spitzer,
\cite{S}, p.~42. The range $\mathcal{R}(  n)$ of the random walk
$S_{ n}$ is the set of sites visited by the walk up to step $n$:
\begin{equation}
\mathcal{R}( n)=\{S_{ 1},\dots,S_{ n} \}.\label{1.1}
\end{equation} As usual, $|\mathcal{R}( n) |$ denotes the
cardinality of the range  up to step $n$.

Dvoretzky and Erdos, \cite{DEr}, show that
\begin{equation}
    \lim_{ n\rar\ff} \log n { |\mathcal{R}(n)|\over  n}=2\pi,\hspace{
.2in}a.s.\label{rp7.2}
\end{equation}
     Le Gall \cite{LeGall6} has obtained a central limit theorem  for
the second order fluctuations of $ |\mathcal{R}(n)|$:
\begin{equation}
     (\log n)^{ 2}
\({ |\mathcal{R}(n)|-E( |\mathcal{R}(n)| )\over
n}\)\stackrel{d}{\rightarrow}-(2\pi)^{ 2}
\ga_{2}( 1)\label{rp7.3}
\end{equation} where $\stackrel{d}{\rightarrow}$ denotes
convergence in law and $\ga_{2}( t)$ is the second order
renormalized self-intersection local time for planar Brownian
motion. See also Le Gall-Rosen \cite{LR}.

In this paper we prove an a.s.\  asymptotic expansion for  $
|\mathcal{R}(n)|$ to any order of accuracy. In order to state our
result we first introduce some notation.
  If $\{W_{ t}\,;\,t\geq 0\}$ is a planar Brownian motion, we define
the
$j$'th  order renormalized intersection local time for $\{W_{
t}\,;\,t\geq 0\}$ as follows.
${\ga}_{ 1}( t)=t$, $\alpha_{1,\ep}(t)=t$  and for $k\geq 2$
\begin{equation}\label{al1e}
\alpha_{k,\ep}(t)=\int_{ 0\leq t_{ 1}\leq\cdots\leq t_{ k}<t}
\prod_{ i=2}^{k} p_{ \ep}(W_{ t_{ i}}-W_{ t_{ i-1}})dt_{ 1}\cdots
dt_{k},
\end{equation}
\begin{equation}\label{2.4a}
\ga_{k}(t)=\lim_{ \ep\rar 0}\sum_{ l=1}^{ k}{k-1\choose l-1}
\(-u_{ \ep}\)^{ k-l} \alpha_{l,\ep}(t),
\end{equation} where
$p_{ t}(x)$ is the density for $W_{ t}$
     and
\[u_{ \ep}=\int_{ 0}^{ \ff} e^{ -t}p_{t+ \ep}(0) \,dt.
\]

Renormalized self-intersection local time was originally  studied by
Varadhan
\cite{Varadhan} for its role in quantum field theory. In  Rosen
\cite{RosenDM} we show that
$\ga_{k}( t)$ can be characterized as the  continuous process of
zero quadratic variation in the decomposition of a natural
Dirichlet process. For further work on renormalized 
self-intersection local times see Dynkin
\cite{Dynkin88}, Le Gall \cite{LeGall-St.Flour}, Bass and
Khoshnevisan
\cite{BK}, Rosen
\cite{RosenJC} and Marcus and Rosen \cite{Marcus-Rosen12}.

To motivate our result define the Wiener sausage of radius $\ep$
as
\begin{equation}
\mathcal{W}_{ \ep}( 0,t)=\{x\in R^{ 2}\,|\, \inf_{ 0\leq s\leq t}|x-
W_{ s}|\leq \ep\}.\label{ws.1}
\end{equation} Letting $m(\mathcal{W}_{ \ep}( 0,t) )$ denote the
area of the Wiener sausage of radius $\ep$, Le Gall
\cite{LG90} shows that for each
$k\geq 1$
\[ (\log n)^{k}\left[ m(\mathcal{W}_{n^{ -1/2}}( 0,1)
)+\sum_{j=1}^{k} ( -1)^{ j}(\textstyle{{1 \over 2\pi}}\log 
n+c)^{-j} \ga_{ j}(1)\right]\to 0,
\hspace{ .2in} a.s.
\] as $n\rar\ff$ where $c$ is a finite constant. Using the heuristic
which associates
$\{S_{[nt]}/\sqrt{n}\,;0\leq t\leq 1 \}\subseteq n^{-1/2}Z^{ 2}
\subseteq R^{ 2}$ with the Brownian motion
$\{W_{ t}\,;\,0\leq t\leq 1\}$, one would expect (note that space
is scaled by
$n^{ -1/2}$) that
$\frac{1}{n} |\mathcal{R}(n)|$ will be `close' to
$m(\mathcal{W}_{n^{-1/2}}( 0,1) )$.

Our main result is the following theorem.

\bt\label{theo-fixasip} Let $S_{ n}=X_{ 1}+\cdots+X_{ n}$ be a
symmetric, strongly aperiodic random walk in $Z^2$ with
covariance matrix equal to the identity and with
$2+\de$ moments for some $\de>0$. On a suitable probability
space we can construct $\{S_{ n}\,;\,n\geq 1\}$ and a planar
Brownian motion
$\{W_{ t}\,;\,t\geq 0\}$ such that for each
$k\geq 1$
\begin{equation} \hspace{ .3in} (\log n)^{k}\left[\frac{1}{n}
|\mathcal{R}(n)|+\sum_{j=1}^{k} ( -1)^{ j}(\textstyle{{1 \over
2\pi}}\log  n+c_{X})^{-j} \ga_{ j}(1,W^{(n)})\right]\to 0,
\hspace{ .2in} a.s.\label{b.10}
\end{equation} where the random variables $\ga_{ 1}(1,W^{(n)}
),\ga_{ 2}(1,W^{(n)}  ),\ldots$ are the renormalized
self-intersection local times (\ref{2.4a}) with $t=1$  for the
Brownian motion
$\{W^{(n)}_t=W_{nt}/\sqrt{n}\,;\,t\geq 0\}$,
\begin{equation} c_X={1 \over 2\pi}\log (\pi^{ 2}/2 )+ {1 \over
(2\pi)^{ 2}}\int_{[-\pi,\pi]^2} {\phi(p)-1+|p|^{ 2}/2\over
(1-\phi(p))|p|^{ 2}/2}\,dp\label{1.5cx}
\end{equation} is a finite constant, and  $\phi(p)=E(e^{ip\cdot
X_1})$
  denotes the characteristic function of $X_1$.
\et

Note that the presence of the constant $c_X$ shows that the
heuristic mentioned before the statement of Theorem
\ref{theo-fixasip} does not completely capture the fine structure
of $|\mathcal{R}(n)|$. (This can already be observed on the level
of (\ref{rp7.3}), see \cite[( 6.r)]{LR}).

We begin our proof in Section \ref{sec-rwilt} where we introduce
renormalized intersection local times $\Ga_{ k,\la}( n)$ for our
random walk.  Let $\ze$ be an independent exponential random
variable of mean 1,  and set $\ze_{
\la}=n$ when
$(n-1 )\la<\ze\leq \la n$. Letting $|\mathcal{R}( \ze_{ \la})|$
denote the cardinality of the range of our random walk killed at
step
$\ze_{\la}$, we derive an $L^{ 2}$ asymptotic expansion for
  $|\mathcal{R}( \ze_{ \la})|$ in terms of the
  $\Ga_{ k,\la}(\ze_{\la})$ as $\la\rar 0$. In Sections
\ref{sec-strongapp}-\ref{sec-approxilt}, on a suitable probability
space,  we construct
  $\{S_{ n}\,;\,n\geq 1\}$ and a planar Brownian motion
$\{W_{ t}\,;\,t\geq 0\}$ and show that in the above $L^{ 2}$
asymptotic expansion for
  $|\mathcal{R}( \ze_{ \la})|$ we can replace  $\la\Ga_{
k,\la}(\ze_{\la})$ by
$\ga_{ k}(\ze,W^{(\la^{ -1})})$, the renormalized intersection
local times  for the planar Brownian motion
$\{W^{(\la^{ -1})}_t=W_{\la^{ -1}t}/\sqrt{\la^{ -1}}\,;\,t\geq
0\}$. After some preliminaries on renormalized intersection local
times for Brownian motion in Section
\ref{sec-brilt}, we show in Section
\ref{sec-rabilt} how our $L^{ 2}$ asymptotic expansion for
  $|\mathcal{R}( \ze_{ \la})|$ leads to an a.s. asymptotic expansion.
The proof of Theorem \ref{theo-fixasip} is completed in Section
\ref{sec-nonrandom} by showing how to replace the random time
$\ze_{ \la}$ by fixed time. Appendix \ref{sec-estimates} derives
some estimates used in this paper.
  Our  methods obviously owe a great deal to Le Gall
\cite{LG90}.

We would like to thank Uwe Einmahl and  David Mason for their
help in Section \ref{sec-strongapp} which describes strong
approximations
  in $L^{ 2}$.

\section{Range and random walk intersection local
times}\label{sec-rwilt}

We first define the non-renormalized random walk intersection
local times  for
$k\geq 2$ by
\bea I_{ k}( n)&=&\sum_{0\leq  i_{ 1}\leq\ldots\leq  i_{ k}< n}
\de(S_{ i_{ 1}}, S_{ i_{ 2}})\cdots \de(S_{ i_{k- 1}}, S_{ i_{ k}})
\label{1.3}\\ &=&\sum_{ x \in Z^{2}}\sum_{0\leq  i_{
1}\leq\ldots\leq i_{ k}< n}
\prod_{ j=1}^{ k}\de(S_{ i_{ j}}, x)\nn
\eea where
\[
\de(i,j)=\left\{ \begin{array}{ll} 1 &\mbox{if $i=j$}\\ 0 &
\mbox{otherwise}
\end{array}
\right.
\] is the usual Kronecker delta function. We set $I_{ 1}( n)=n$ so
that also
$I_{ 1}( n)=\sum_{ x \in Z^{2}}\sum_{0\leq  i< n}
\de(S_{ i}, x)$. (One might also take as a definition of the 
intersection local time the quantity
$\sum_{0<  i_{ 1}<\ldots<  i_{ k}< n}
\de(S_{ i_{ 1}}, S_{ i_{ 2}})\cdots \de(S_{ i_{k- 1}}, S_{ i_{ k}}).$ 
The definition in (\ref{1.3})  is more convenient for our purposes, 
and we see by (\ref{1.6}) that either definition leads to the same
value for
$\Gamma_{k,\lambda}(n)$.)

Let $q_{ n}( x)$ be the transition function for
$S_{ n}$ and  let
\begin{equation} G_{ \la}( x)=\sum_{ j=0}^{ \ff}e^{ -j\la}q_{ j}(
x).\label{1.4}
\end{equation} We will show in Lemma \ref{lem-green} below that
\begin{equation} g_{ \la}:=G_{ \la}( 0)= {1 \over 2\pi}\log (
1/\la)+c_{X}+O( \la^{ \de }\log (1/\la))
\hspace{.1in} \mbox{ as    }
\hspace{ .1in} \la\rar 0,\label{1.5}
\end{equation} where $c_X$ is defined in (\ref{1.5cx}). We show
in (\ref{45.31}) that for any
$q>1$
\begin{equation}
\sum_{ x \in Z^{2}}(G_{ \la}( x))^{ q}=O(\la^{ -1} )\hspace{ .1in}
\mbox{ as    }
\hspace{ .1in} \la\rar 0.\label{1.5b}
\end{equation} Note also that
\begin{equation}\label{2.2AA}
\sum_{x\in Z^2} G_\lambda(x)=\sum_{j=0}^\infty
e^{-j\lambda}=\frac{1}{1-e^{-\lambda}}.
\end{equation} We now define the renormalized random walk
intersection local times
\begin{eqnarray}&&
\Ga_{ k,\la}( n)\nn\\ && =\sum_{0\leq  i_{ 1}\leq\ldots\leq  i_{
k}< n}
\{ \de(S_{ i_{ 1}}, S_{ i_{ 2}})-g_{ \la}\de( i_{ 1},  i_{ 2})\}
\cdots
\{ \de(S_{ i_{k- 1}}, S_{ i_{ k}})-g_{ \la}\de( i_{k- 1},  i_{ k})\}
\nonumber\\ &&\qquad =\sum_{j=1}^{ k}{k-1\choose j-1}( -1)^{
k-j}g^{ k-j}_{ \la}I_{ j}( n).
\label{1.6}
\end{eqnarray}

Let $\ze$ be an independent exponential random variable of
mean 1, and set $\ze_{\la}=n$ when
$(n-1 )\la<\ze\leq \la n$.
     $\ze_{ \la}$ is then a geometric random variable with
$P( \ze_{ \la}>n)=e^{ -\la n}$.  Note that $\zeta_{1/j}=n$ if
$(n-1)/j<\zeta\leq n/j.$ By $\mathcal{R}( \ze_{ \la})$ we mean 
the range of our random walk killed at step $\ze_{\la}$.

    In this section we prove the following lemma.
\begin{lemma}\label{lem-rwasip} For each
$k\geq 1$
\begin{equation}
\lim_{  \la\to 0}\la g^{ k}_{ \la}
\(|\mathcal{R}( \ze_{  \la})|-\sum_{ j=1}^{ k}( -1)^{ j-1}g^{ -j}_{
\la}\Ga_{ j,\la}( \ze_{  \la})\)=0\hspace{ .1in}\mbox{ in
}L^2.\label{1.7}
\end{equation}
\end{lemma}

{\bf Proof of Lemma \ref{lem-rwasip}:}  Define
\[ T_x=\min\{n\geq 0: S_n=x\},
\] the first hitting time to $x$. We will use the fact that
\begin{equation} P(T_{x}<\ze _{\la})={G_{ \la}( x) \over G_{ \la}(
0)},\label{1.11}
\end{equation} which follows from the strong Markov property:
\begin{eqnarray} G_{ \la}( x)&=&E\(\sum_{ j=0}^{ \ff}1_{ x}( S_{
j})1_{
\{j<\ze _{\la}
\}}\)\label{}\\ &= &E\(\sum_{ j=T_{x}}^{ \ff}1_{ x}( S_{ j})1_{
\{j<\ze _{\la}
\}}\)
\nonumber\\ &= &E\(1_{ \{T_{x}<\ze _{\la} \}},\,\lc\sum_{
j=0}^{
\ff}1_{ x}( S_{  j})1_{ \{j<\ze _{\la}
\}}\rc\circ\th_{T_{x} }\)
\nonumber\\ &= &P(T_{x}<\ze _{\la})G_{ \la}( 0).
\nonumber
\end{eqnarray}

To prove our lemma we square the expression inside the
parentheses in (\ref{1.7}) and then take expectations. We first
show that
\begin{equation}\qquad E\(|\mathcal{R}( \ze_{  \la})|^{
2}\)=2\sum_{ j=2}^{ 2k}( -1)^{ j}g^{  -j}_{ \la}
\sum_{ x,y\in Z^{2}} G_{ \la}( x)(G_{ \la}( x-y))^{ j-1} +O(\la^{
-2}g^{ -(2k+1)}_{ \la}).\label{1.8}
\end{equation} To this end we first note that
\begin{equation} |\mathcal{R}( \ze_{  \la})|=\sum_{ x\in
Z^{2}}1_{
\{ T_{ x}<\ze _{
\la}\}}\label{1.9}
\end{equation} so that
\bea E\(|\mathcal{R}( \ze_{  \la})|^{ 2}\)&=& \sum_{ x,y\in
Z^{2}}P(T_{ x},T_{ y}<\ze _{
\la})\label{1.10}\\ &=& \sum_{ x \in Z^{2}}P(T_{ x}<\ze
_{\la})+2\sum_{ x\neq y\in Z^{2}}P(T_{ x}<T_{ y}<\ze _{\la}).\nn
\eea Using (\ref{1.11}) we have that
\begin{equation}\qquad
\sum_{ x \in Z^{2}}P(T_{ x}<\ze _{\la})=\sum_{ x \in Z^{2}}{G_{
\la}( x)
\over g_{ \la}} ={1 \over (1-e^{ -\la})g_{ \la}}=O(\la^{ -1}g^{ -1}_{
\la}).\label{1.12}
\end{equation} To evaluate $\sum_{ x\neq y\in Z^{2}}P(T_{
x}<T_{ y}<\ze _{\la})$ we first introduce some notation. For any
$u\neq v\in Z^{2}$ define inductively
\bea && A^{ 1}_{ u,v}=T_{ u}\label{1.13}\\ && A^{ 2}_{ u,v}=A^{
1}_{ u,v}+T_{ v}\circ\th_{ A^{ 1}_{ u,v}}\nn\\ && A^{ 3}_{
u,v}=A^{ 2}_{ u,v}+T_{ u}\circ\th_{ A^{ 2}_{ u,v}}\nn\\ && A^{
2k}_{ u,v}=A^{ 2k-1}_{ u,v}+T_{ v}\circ\th_{ A^{ 2k-1}_{
u,v}}\nn\\ && A^{ 2k+1}_{ u,v}=A^{ 2k}_{ u,v}+T_{ u}\circ\th_{
A^{ 2k}_{ u,v}}.\nn
\eea

We observe that for any $x\neq y$
\begin{eqnarray} P^{ z}(T_{x}<T_{ y}<\ze _{\la})&=&P^{ z}(A^{
1}_{ x,y}<A^{ 2}_{  x,y}<\ze _{\la}) -P^{ z}(T_{ y}<A^{ 1}_{
x,y}<A^{ 2}_{ x,y}<\ze _{\la})
\nn\\ &=& P^{ z}(A^{ 2}_{ x,y}<\ze _{\la}) -P^{ z}(T_{ y}<A^{ 1}_{
x,y}<A^{ 2}_{ x,y}<\ze _{\la})\label{1.18}
\end{eqnarray} and
\begin{eqnarray} && P^{ z}(T_{ y}<A^{ 1}_{ x,y}<A^{ 2}_{ x,y}<\ze
_{\la})\label{1.19}\\ &&=P^{ z}(A^{ 1}_{ y,x}<A^{ 2}_{ y,x}<A^{
3}_{ y,x}<\ze _{\la}) -P^{ z}(T_{ x}< A^{ 1}_{ y,x}<A^{ 2}_{
y,x}<A^{ 3}_{ y,x}<\ze _{\la})
\nn\\&& = P^{ z}(A^{ 3}_{ y,x}<\ze _{\la}) -P^{ z}(T_{ x}<A^{ 1}_{
y,x};\,A^{ 3}_{ y,x}<\ze _{\la}).\nn
\end{eqnarray} Proceeding inductively we find that
\begin{eqnarray} P^{ z}(T_{x}<T_{ y}<\ze _{\la})&=&\sum_{
j=1}^{ k}P^{ z}(A^{ 2j}_{  x,y}<\ze _{\la}) -\sum_{ j=1}^{
k}P^{z}(A^{ 2j+1}_{ y,x}<\ze _{\la})
\nn\\ && +P^{ z}(T_{ x}<A^{ 1}_{ y,x};\,A^{ 2k+1}_{ y,x}<\ze
_{\la}).\label{1.20}
\end{eqnarray} Using (\ref{1.11}) and the strong Markov
property we see that
\begin{eqnarray} P(T_{x}<T_{ y}<\ze _{\la})&=&\sum_{ j=1}^{ k}
g^{ -2j}_{\la}G_{ \la}( x) (G_{ \la}( y-x))^{ 2j-1}\nn\\ && -\sum_{
j=1}^{ k}g^{ -(2j+1)}_{\la}G_{ \la}( y) (G_{ \la}( x-y))^{ 2j}
\nn\\ && +P(T_{ x}<A^{ 1}_{ y,x};\,A^{ 2k+1}_{ y,x}<\ze _{\la})
\label{1.21}
\end{eqnarray} and that
\bea && P(T_{ x}<A^{ 1}_{ y,x};\,A^{ 2k+1}_{ y,x}<\ze
_{\la})\label{1.22}\\ &&\leq P(\,A^{ 2k+2}_{ x,y}<\ze _{\la})=g^{
-(2k+2)}_{\la}G_{ \la}(  x) (G_{
\la}( y-x))^{ 2k+1}.\nn
\eea (\ref{1.8}) then follows using (\ref{1.5}) and (\ref{1.5b}).

We next observe that
\begin{eqnarray} && E( I_{ n}( \ze _{\la})I_{ m}( \ze
_{\la}))\label{1.23x}\\ && =\sum_{ x,y\in Z^{2}}E\(\sum_{0\leq 
i_{ 1}\leq\ldots\leq  i_{ n}<\ze _{\la}}
\prod_{ j=1}^{ n}\de(S_{ i_{ j}}, x)\sum_{0\leq  l_{
1}\leq\ldots\leq l_{ m}<\ze _{\la}}
\prod_{ k=1}^{ m}\de(S_{ l_{ k}}, x)\)\nonumber\\ && =\sum_{
x,y\in Z^{2}}\sum_{ \pi}E\(\sum_{0\leq  i_{ 1}\leq\ldots\leq  i_{
n+m}<\ze _{\la}}
\prod_{ j=1}^{ n+m}\de(S_{ i_{ j}}, \pi( j))\)\nonumber
\end{eqnarray} where the inner sum runs over all maps $\pi: \,\{
1,2,\ldots,n+m\}\mapsto \{ x,y\}$ such that $|\pi^{ -1}(
x)|=m,\,|\pi^{ -1}( y)|=n$. Thus
\begin{eqnarray} && E( I_{ n}( \ze _{\la})I_{ m}( \ze
_{\la}))\label{1.23}\\ && =\sum_{ x,y\in Z^{2}}\sum_{
\pi}\sum_{0\leq  i_{ 1}\leq\ldots\leq  i_{ n+m}<\ff}E\(
\prod_{ j=1}^{ n+m}\de(S_{ i_{ j}}, \pi( j))\)e^{ -\la i_{
n+m}}\nonumber\\ && =\sum_{ x,y\in Z^{2}}\sum_{ \pi}\prod_{
j=1}^{n+m } G_{ \la}( \pi( j)-\pi( j-1))\nonumber
\end{eqnarray} where $\pi( 0)=0$. By (\ref{1.5}) and (\ref{1.5b})
the sum over $x=y$ is $O(\la^{ -1}g^{n+m}_{ \la})$.
     Consider then $x\neq y$. When we look at the definition
(\ref{1.6}) of
$\Ga_{ k,\la}( n)$ we see that the effect of replacing
$ I_{ n}( \ze _{\la})I_{ m}(
\ze _{\la})$ in (\ref{1.23}) by $\Ga_{ n,\la}( \ze _{\la})\Ga_{ 
m,\la}(
\ze _{\la})$ is to eliminate all maps $\pi$ in which $\pi( j)=\pi( 
j-1)$ for some
$j$.  Thus, up to an error which is
$O(\la^{ -1}g^{n+m}_{ \la})$, (which comes from $x=y$), we have
\begin{equation}\qquad E(\Ga_{ n,\la}( \ze _{\la})\Ga_{ m,\la}(
\ze _{\la}))=\left\{ \begin{array}{ll} 2\sum_{ x,y\in Z^{2}}G_{
\la}( x) (G_{
\la}( y-x))^{2n-1} &\mbox{if $m=n$}\\
\sum_{ x,y\in Z^{2}}G_{ \la}( x) (G_{ \la}( y-x))^{2n-1\pm 1}
&\mbox{if
$m=n\pm 1$}\\ 0 & \mbox{otherwise.}
\end{array}
\right.\label{1.24}
\end{equation} Consequently  up to errors which are
$O(\la^{ -1}g^{2k}_{ \la})$
\begin{eqnarray} && E\(\lc\sum_{ j=1}^{ k}( -1)^{ j-1}g^{ -j}_{
\la}\Ga_{ j,\la}( \ze_{  \la})\rc^{ 2}\)\label{1.25}\\ && =\sum_{
n,m=1}^{ k}( -1)^{n+m}g^{ -( n+m)}_{ \la} E(\Ga_{ n,\la}( \ze
_{\la})\Ga_{ m,\la}( \ze _{\la}))\nonumber\\ && =2\sum_{
n=1}^{ k}g^{ -2n}_{ \la}
\sum_{ x,y\in Z^{2}}G_{ \la}( x) (G_{ \la}(
y-x))^{2n-1}\nonumber\\ && ~~~ -2\sum_{ n=2}^{ k}g^{
-(2n-1)}_{ \la}
\sum_{ x,y\in Z^{2}}G_{ \la}( x) (G_{ \la}(
y-x))^{2n-2}\nonumber\\ && =2\sum_{ j=2}^{ 2k}( -1)^{ j}g^{
-j}_{ \la}
\sum_{ x,y\in Z^{2}} G_{ \la}( x)(G_{ \la}( x-y))^{ j-1}\nn
\end{eqnarray}

To handle the cross product terms we define the random measure
on
$Z^{ n}_{ +}$
\bea
\La_{ n,y}( B)&=&\sum_{\{ 0\leq  i_{ 1}\leq\ldots\leq  i_{
n}<\ze_{
\la}\}\cap B}
\prod_{ j=1}^{ n}\de(S_{ i_{ j}}, y).
\label{1.26}
\eea Using the notation $i_{ 0}=0, i_{ n+1}=\ze_{  \la}$ we have
\begin{eqnarray} E(|\mathcal{R}( \ze_{  \la})|I_{ n}( \ze _{\la}))
&=&\sum_{x,y\in Z^2} \sum_{0\leq i_1\leq \cdots \leq i_n\leq
\ze_\la} 1_{(T_x<\zeta_\la)} \prod_{j=1}^n \de(S_{i_j},y) \nn\\
&=&
\sum_{ x,y\in Z^{2}}\sum_{ j=0}^{ n}E(\La_{ n,y}( \{i_{ j}\leq T_{
x}< i_{ j+1}
\})).\label{1.27}
\end{eqnarray} As above we have that
\begin{eqnarray} &&
\La_{ n,y}( \{i_{ j}\leq T_{ x}< i_{ j+1}\})
\label{1.28}\\ && =\La_{ n,y}( \{i_{ j}+T_{ x}\circ\th_{ i_{ j}}< i_{
j+1}\})\nonumber\\ && -\sum_{ l=0}^{ j-1}
\La_{ n,y}( \{i_{ l}\leq T_{ x}< i_{ l+1};\,i_{ j}+T_{ x}\circ\th_{ 
i_{ j}}< i_{ j+1}\})
\nonumber
\end{eqnarray} and inductively we find that
\begin{eqnarray} &&
\sum_{ x\neq y\in Z^{2}}\sum_{ j=0}^{ n}E(\La_{ n,y}( \{i_{
j}\leq  T_{ x}< i_{ j+1}
\}))
\label{1.29}\\ && =\sum_{ x\neq y\in Z^{2}}\sum_{ m=1}^{
n+1}( -1)^{ m-1}
\sum_{ |A|=m}E\(\La_{ n,y}( \cap_{ j\in A}\{i_{ j}+T_{
x}\circ\th_{ i_{ j}} < i_{j+1}\})\)
\nonumber
\end{eqnarray} where the inner sum runs over all nonempty
$A\subseteq \{ 0,1,\ldots,n\}$. Using (\ref{1.11}) and the Markov
property we see that
\begin{eqnarray} &&
\sum_{ x\neq y\in Z^{2}}\sum_{ j=0}^{ n}E(\La_{ n,y}( \{i_{
j}\leq  T_{ x}< i_{ j+1}
\}))
\label{1.30h}\\ && =\sum_{ x\neq y\in Z^{2}}\sum_{ m=1}^{
n+1}( -1)^{ m-1}
\sum_{ |A|=m}g^{ -m}_{ \la}\prod_{ j=1}^{n+m } G_{ \la}( \si_{
A}( j)-\si_{ A}( j-1))
\nonumber
\end{eqnarray} where $\si_{ A}(0)=0$ and $\si_{ A}( j)$ is the
$j$'th element in the  ordered set obtained by taking
$n$
     $y$'s and inserting, for each $l\in A$, an $x$ between the
$l$'th  and
$(l+1)$'st $y$. Estimating the contribution from $x=y$ we find
that
\begin{eqnarray} && E(|\mathcal{R}( \ze_{  \la})|I_{ n}( \ze
_{\la}))
\label{1.30}\\ && =\sum_{ x,y\in Z^{2}}\sum_{ m=1}^{ n+1}(
-1)^{ m-1}
\sum_{ |A|=m}g^{ -m}_{ \la}\prod_{ j=1}^{n+m } G_{ \la}( \si_{
A}( j)-\si_{ A}( j-1))\nonumber\\ &&\hspace{ .5in}+O(\la^{ -2}g^{
-(2k+1)}_{ \la}).
\nonumber
\end{eqnarray} Once again we see that the effect of replacing $
I_{ n}(
\ze _{\la})$ in (\ref{1.30}) by $\Ga_{ n,\la}( \ze _{\la})$ is to
eliminate all sets
$A$ such that $\si_{ A}( j)=\si_{ A}( j-1)$ for some $j$. Thus we
have
\begin{eqnarray} && E(|\mathcal{R}( \ze_{  \la})|\Ga_{ n,\la}( \ze
_{\la}))
\label{1.31}\\ && =2( -1)^{ n-1}\sum_{ x,y\in Z^{2}} g^{ -n}_{
\la}G_{
\la}( x) (G_{ \la}( x-y))^{ 2n-1}
\nonumber\\ && +( -1)^{ n}\sum_{ x,y\in Z^{2}}g^{ -(n-1)}_{
\la}G_{
\la}( x) (G_{ \la}( x-y))^{ 2n-2}
\nonumber\\ && +( -1)^{ n}\sum_{ x,y\in Z^{2}}g^{ -(n+1)}_{
\la}G_{
\la}( x) (G_{ \la}( x-y))^{ 2n}\nn\\ &&\hspace{ .5in}+O(\la^{
-2}g^{ -(2k+1)}_{ \la}).
\nonumber
\end{eqnarray} Consequently
\begin{eqnarray} && E\(|\mathcal{R}( \ze_{  \la})|\sum_{ n=1}^{
k}( -1)^{ n-1}g^{ -n}_{ \la}\Ga_{ n,\la}( \ze_{ 
\la})\)\label{1.32}\\ && =2\sum_{ n=1}^{ k}g^{ -2n}_{ \la}
\sum_{ x,y\in Z^{2}}G_{ \la}( x) (G_{ \la}(
y-x))^{2n-1}\nonumber\\ && -\sum_{ n=2}^{ k}g^{ -(2n-1)}_{ \la}
\sum_{ x,y\in Z^{2}}G_{ \la}( x) (G_{ \la}(
y-x))^{2n-2}\nonumber\\ && -\sum_{ n=1}^{ k}g^{ -(2n+1)}_{
\la}
\sum_{ x,y\in Z^{2}}G_{ \la}( x) (G_{ \la}( y-x))^{2n}\nonumber\\
&&\hspace{ .5in}+O(\la^{ -2}g^{ -(2k+1)}_{ \la})
\nonumber\\ && =2\sum_{ j=2}^{ 2k}( -1)^{ j}g^{ -j}_{ \la}
\sum_{ x,y\in Z^{2}} G_{ \la}( x)(G_{ \la}( x-y))^{ j-1}\nn\\
&&\hspace{ .5in}+O(\la^{ -2}g^{ -(2k+1)}_{ \la}).
\nonumber
\end{eqnarray} Our lemma then follows from (\ref{1.8}),
(\ref{1.25}) and (\ref{1.32}).
\qed

\section{Strong approximation in $L^{ 2}$}\label{sec-strongapp}

\begin{lemma}\label{lem-einmas} Let $X$ be a $R^{2}$ valued
random vector with mean zero and covariance matrix equal to the
identity $I.$ Assume that for some $2<p<4$,
$E|X|^{p}<\infty $. Given $n\geq 1$  one can construct on a
suitable probability space two sequences of independent random
vectors
$X_{1},\dots ,X_{n}$ and $Y_{1},\dots ,Y_{n},$ where each
$X_{i}\stackrel{d}{=}X$ and the $Y_{i}$'s  are standard normal
random vectors such that
\[
\Bigl\| \max_{1\leq k\leq n}\Bigl| \sum_{i=1}^{k}\left(
X_{i}-Y_{i}\right)
\Bigr|\, \Bigr\| _{2}=O\left( n^{1/p }\right).
\]
\end{lemma}

{\bf Proof of Lemma \ref{lem-einmas}:}  By equation (3.3) of
\cite{Einmahl} we can find  a constant $c_1$ and such $X_i$ and
$Y_i$ so that
\[ P\left\{\max_{1\leq k\leq n} |\sum_{i=1}^k
(X_i-Y_i)|>x\right\}\leq c_1 n x^{-p} E|X|^p
\]  Write $Z_n$ for $\max_{1\leq k\leq n} |\sum_{i=1}^k
(X_i-Y_i)|$. Since probabilities are bounded by 1, we have for any
$a>0$
\begin{eqnarray}
\| Z_n \|_2^2 &=& \int_0^\infty x P(Z_n>x)\, dx\nn\\ & \leq &
\int_0^a x\, dx +c_1n \int_a^\infty x^{-p+1}\, dx\nn\\ & \leq &
c_2(a^2+n a^{2-p}).\nn
\end{eqnarray} If we set $a=n^{1/p}$ and take square roots of
both sides, we have our result.
\qed

Using the lemma  we can readily construct two i.i.d.\ sequences
$\left\{ X_{i}\right\} _{i\geq 1}$,  and $\left\{ Y_{i}\right\}
_{i\geq 1},$ where the
$X_i$ are equal in law to $X$ and the $Y_i$ are standard normal,
such that for some constant $C>0$ and any $m\geq 0,$
\[
\left\| \max_{2^{m}\leq k<2^{m+1}}\left|
\sum_{i=2^{m}}^{k}\left( X_{i}-Y_{i}\right) \right|\,  \right\|
_{2}\leq C\left( 2^{m}\right) ^{1/p}.
\] We see then that for any $2^{m}\leq \left[ nt\right] <2^{m+1},$
\[
\left\| \sum_{i=1}^{\left[ nt\right] }\left( X_{i}-Y_{i}\right)
\right\| _{2}\leq
\sum_{j=0}^{m}\left\| \max_{2^{m}\leq k<2^{m+1}}\left|
\sum_{i=2^{m}}^{k}\left( X_{i}-Y_{i}\right) \right| \right\| _{2},
\] which for some $D>0$ is
\[
\leq \sum_{j=0}^{m}C\left( 2^{j}\right) ^{1/p }\leq D\left(
nt\right) ^{1/p }.
\] Now choose a Brownian motion $W$ such that for $m\geq 1,$
\[ W\left( m\right) =\sum_{i=1}^{m}Y_{j}.
\] Noting that
\[
\left\| W\left( \left[ mt\right] \right) -W\left( mt\right) \right\|
_{2}\leq
\left\| \sup_{0\leq s\leq 1}\left| W\left( s\right) \right|
\right\| _{2}:=M.
\] we see that for any $t>0$
\begin{eqnarray}
\left\| \frac{S\left( \left[ mt\right] \right) -W\left( mt\right)
}{\sqrt{m}}
\right\| _{2}&\leq& D\left( mt\right) ^{1/p
}m^{-1/2}+Mm^{-1/2}\nn\\ &=&O\left( m^{(1/p)-(1/2) }\left(
{t}^{1/p}+1\right) \right) ,
\label{strongapprox}
\end{eqnarray} where
\[ S\left( \left[ mt\right] \right) =\sum_{i\leq \left[ mt\right]
}X_{i}.
\]

\section{Spatial H\"older continuity for renormalized intersection
local times}\label{sec-asiprilt}

If $\{W_{ t}\,;\,t\geq 0\}$ is a planar Brownian motion,  set
$\ol \al_{1,\ep}( t)=t$ and for $k\geq 2$ and $x=(x_{ 2},\ldots,x_{ k}
)\in (R^{2})^{k-1}$ let
     \begin{equation}
\ol \al_{k,\ep}(t,x)=\int_{ 0\leq t_{ 1}\leq\cdots\leq t_{ k}<t} \prod_{
i=2}^{ k} p_{ \ep}(W_{ t_{ i}}-W_{ t_{ i-1}}-x_{ i})dt_{ 1}\cdots
dt_{ k}.\label{7.3}
\end{equation} When $x_{ i}\neq 0$ for all $i$ and $\ze$ is an
independent exponential random variable with mean 1, the limit
\begin{equation}
\ol \al_{k}(\ze,x)=\lim_{ \ep\rar 0}\,\ol \al_{k,\ep}(\ze,x)\label{7.3l}
\end{equation} exists. When $x_{ i}\neq 0$ for all $i$ set
\begin{equation}
\ol \ga_{k}(\ze,x)=\sum_{A\subseteq \{ 2,\ldots,k\}}( -1)^{|A|}
\(\prod_{ i\in A} u^{ 1}( x_{ i})\)\ol \al_{k-|A|}(\ze,x_{ A^{
c}}),\label{7.5}
\end{equation} where
\begin{equation} u^{ 1}( y)=\int_{ 0}^{ \ff} e^{ -t}p_{t}(y)
\,dt,\label{7.4}
\end{equation}
$p_{ t}(x)$ is the density for $W_{ t}$, and $x_{A^c}=(x_{i_1},
\ldots, x_{i_{k-|A|}})$ with
$i_1<i_2<\cdots <i_{k-|A|}$ and $i_j\in \{2, \ldots, k\}-A$ for each
$j$, that is, the vector $(x_2, \ldots, x_k)$ with all terms that have
indices in
$A$ deleted. In \cite{RosenJC} it is shown that for some $\de>0$
and all
$m$
\begin{equation} E\(\Big|\ol \ga_{k}(\ze,x)-\ol \ga_{k}(\ze,y)\Big|^{
m}\)\leq C|x-y|^{\bar{\de}  m}.\label{7.6}
\end{equation}

As before, set $I_{ 1}( n)=n$ and for $k\geq 2$ and $x=(x_{
2},\ldots,x_{ k} )\in Z^{2}$ let
\be \ol I_{ k}( n,x)=\sum_{0\leq  i_{ 1}\leq\ldots\leq  i_{ k}< n}
\de(S_{ i_{ 2}}- S_{ i_{ 1}}-x_{ 2})\cdots \de(S_{ i_{k}}-S_{ i_{ k-
1}}-x_{ k}).
\label{7.2}
\ee and for $x\in \sqrt{\la}Z^{2}$ let
\be
\ol \Ga_{ k,\la}( n,x)=\sum_{A\subseteq \{ 2,\ldots,k\}}(
-1)^{|A|}\prod_{ i\in A} G_{ \la}( x_{ i}/\sqrt{\la}) \ol I_{ k-|A|}( n,
x_{ A^{ c}}/\sqrt{\la}).
\label{7.1}
\ee Note that $\Ga_{ k,\la}( n)=\ol \Ga_{ k,\la}( n,0)$.

\begin{lemma}\label{lem-hc} For  any $j\geq 1$ we can find some
$\rho,\bar{\de}>0$ such that uniformly  in $\la>0$
\begin{equation}
\sup_{ |y|\leq \la^{ \rho}}E\(\Big|\la\ol \Ga_{j,\la}(
\ze_{\la},y)-\la\Ga_{j,\la}(
\ze_{\la})\Big|^{ 2}\)\leq C\la^{\bar{\de} }.\label{7.7}
\end{equation}
\end{lemma}

{\bf Proof of Lemma \ref{lem-hc}:}  We begin by considering
\begin{equation} E\(\ol \Ga_{k,\la}( \ze_{\la},x^{1})\ol \Ga_{k,\la}(
\ze_{\la},x^{2})\)\label{4h.1}
\end{equation} for $x^{ i}\in ( Z^{2})^{ k-1}.$

If $h$ is a function which depends on the variable $x$, let
\[\mathcal{D}_x h=h( x)-h(0).\] Let $\mathcal{S}$ be the set of all
maps
$s:\,\{  1,2,\ldots,2k\}\mapsto \{ 1,2\}$ with
$|s^{ -1}( j)|=k,\,1\leq j\leq 2$, and let
     $B_{ s}=\{ i\,|\,s( i)=s( i-1)\}$ and $c( i)=|\{j\leq  i\,|\,s(j)=s( i)
\}|$.

Using the Markov property as in Lemma 5 of \cite{RosenJC} we
can then show that
\begin{eqnarray} && E\(\ol \Ga_{k,\la}( \ze_{\la},x^{1})\ol \Ga_{k,\la}(
\ze_{\la},x^{2})\)\label{4h.2}\\ && =\sum_{s\in\mathcal{S}}\(
\prod_{ i\in B_{ s}}G_{\la}( x^{ s( i)}_{  c( i)}/\sqrt{\la})\)
     \sum_{ \stackrel{z_{ i}\in Z^{ 2}}{ i=1,2}}\( \prod_{ i\in
B_{s}}\mathcal{D}_{ x^{ s( i)}_{ c(  i)}/\sqrt{\la}}\)\nonumber\\
&&
\prod_{ i\in B^{ c}_{ s}}G_{\la}\(z_{s( i) }+\sum_{ j=2}^{ c( i)}
x^{  s( i)}_{ j}/\sqrt{\la} -(z_{s( i-1) }+\sum_{ j=2}^{ c( i-1)} x^{ s(
i-1)}_{ j}/\sqrt{\la} )\).\nn
\end{eqnarray}

     Fix $s\in\mathcal{S}$ and note then that the corresponding
summand will be $0$ unless
$x^{ s( i)}_{ c( i)}\neq 0$ for all $i\in B_{s}$. Note that by 
definition of
$B_{s}^{ c}$ we necessarily have that the last line in (\ref{4h.2}) 
is of the form
\begin{equation} G_{\la}\(z_{1}\)\prod_{ i\in B^{ c}_{ s},i\neq
1}G_{\la}\(z_{1} -z_{2}+a_{ i} \)\label{4h.2m}
\end{equation} where the $a_{ i} $ are linear combinations of
$x^{1},x^{2}$ but do not involve
$z_{1} ,z_{2}$. Then we observe that the effect of applying  each
$\mathcal{D}_{ x^{ s( i)}_{ c( i)}/\sqrt{\la}}$ to the product on
the  last line of (\ref{4h.2}) is to generate a sum of several terms
in each of which  we have one factor of the form
$\mathcal{D}_{ x^{ s( i)}_{ c( i)}/\sqrt{\la}} G_{ \la}$. Thus
schematically we can write the contribution of such a term as
\begin{eqnarray} &&
   \( \prod_{ i\in B_{ s}}G_{\la}( x^{ s( i)}_{  c(i)}/\sqrt{\la})\)
     \sum_{ \stackrel{z_{ i}\in Z^{ 2}}{ i=1,2}}
G_{\la}\(z_{1}\)\prod_{ i\in B^{ c}_{ s},i\neq
1}\Delta_{ A_{ i}}G_{\la}\(z_{1} -z_{2}+a_{ i} \)\label{4h.2d}
\end{eqnarray} where each $\Delta_{ A_{ i}}$ is a product of $k_{
i}$ difference operators of the form
$\Delta_{x^{j}_{l}/\sqrt{\la}}$, and we have
$\sum_{ i\in B^{ c}_{ s}}k_{ i}=|  B_{ s} |$.  If $B_{ s}\neq
\emptyset
$ and if there is only one term in the last product on the right of
(\ref{4h.2d}), it is easily seen that the sum over $z_{2}$ gives $0$.
Thus the product contains at least $2$ terms and then by
   Lemma
\ref{lem-momgreen} we can see that for some $C<\ff$ and
$\nu>0$ independent of everything
\begin{eqnarray} &&
     \Bigg|    \sum_{ \stackrel{z_{ i}\in Z^{ 2}}{ i=1,2}}
G_{\la}\(z_{1}\)\prod_{ i\in B^{ c}_{ s},i\neq 1}\Delta_{ A_{
i}}G_{\la}\(z_{1} -z_{2}+a_{ i} \)\Bigg|
    \leq C\la^{ -2}
     \,\prod_{ i\in B_{ s}}|x^{ s( i)}_{ c( i)} |^{\nu}.\label{4h.3}
\end{eqnarray}

With these results, we now turn to the bound (\ref{4h.1}). For ease
of exposition we use $y^{i}$ to denote the $y$ in the $i$'th factor;
in the end we will set $y^{i}=y$.  For ease of exposition we assume
that
$y$  differs  from $0$ only in the
$v$'th coordinate, and we set $a=y_{v}$. (The general case is then
easily handled).

We again use Lemma 5 of
\cite{RosenJC} to expand
\begin{equation} E\(\(\ol \Ga_{k,\la}( \ze_{\la},y^{ 1})-\Ga_{k,\la}(
\ze_{\la})\)\(\ol \Ga_{k,\la}( \ze_{\la},y^{ 2})-\Ga_{k,\la}(
\ze_{\la})\)\)\label{4h.15}
\end{equation}
     as a sum of many terms of the form
\begin{eqnarray} &&
\sum_{s\in\mathcal{S}}\(\prod_{i=1}^{2}
\DD_{y^{i}_{v}/\sqrt{\la}}\)
\( \prod_{ i\in B_{s}}G_{\la}( x^{ s( i)}_{ c( i)}/\sqrt{\la})\)
     \sum_{ \stackrel{z_{ i}\in Z^{ 2}}{1=1,2}}\( \prod_{ i\in
B_{s}}\mathcal{D}_{ x^{ s( i)}_{ c(  i)}/\sqrt{\la}}\)\label{4h.4}\\
&&
\prod_{ i\in B^{ c}_{ s}}G_{\la}\(z_{s( i) }+\sum_{ j=2}^{ c( i)}
x^{  s( i)}_{ j}/\sqrt{\la} -(z_{s( i-1) }+\sum_{ j=2}^{ c( i-1)} x^{ s(
i-1)}_{ j}/\sqrt{\la} )\).\nn
\end{eqnarray}
     where now $x^{i}$ is variously $y^{i}$ or
$0$. For fixed $s\in\mathcal{S}$ we can expand the corresponding
term as a sum of terms of the form
\begin{eqnarray} &&
\lc\(\prod_{k\in F} \DD_{y^{k}_{v}/\sqrt{\la}}\)
\( \prod_{ i\in B_{s}}G_{\la}( x^{ s( i)}_{ c(
i)}/\sqrt{\la})\)\rc\label{4h.5}\\ &&
     \sum_{ \stackrel{z_{ i}\in Z^{ 2}}{i=1,2}}
\(\prod_{k\in F^{ c}} \DD_{y^{k}_{v}/\sqrt{\la}}\)
\( \prod_{ i\in B_{s}}\mathcal{D}_{ x^{ s( i)}_{ c(
i)}/\sqrt{\la}}\)\nn\\ &&
\prod_{ i\in B^{ c}_{ s}}G_{\la}\(z_{s( i) }+\sum_{ j=2}^{ c( i)}
x^{  s( i)}_{ j}/\sqrt{\la} -(z_{s( i-1) }+\sum_{ j=2}^{ c( i-1)} x^{ s(
i-1)}_{ j}/\sqrt{\la} )\).\nn
\end{eqnarray} where $F$ runs through the subsets of $\{
1,\ldots,n\}$. Note that  the first line will be $0$ unless for each
$k\in F$ we have that
$y^{k}_{v}= x^{ s( i)}_{  c( i)}$ for some
$i\in B_{s}$. In particular
\begin{equation} |F|\leq |B_{s}|.\label{4h.5a}
\end{equation}Using the fact that
\begin{equation} G_{\la}( x)\leq c\log ( 1/\la)\label{4h.6}
\end{equation} we can bound the first line of (\ref{4h.5}) by
$(c\log ( 1/\la))^{ |B_{ s}|}$. As before, see in particular
(\ref{4h.3}), we  can obtain the bound
\begin{eqnarray} &&
\Bigg| \sum_{ \stackrel{z_{ i}\in Z^{ 2}}{i=1,2}}
\(\prod_{k\in F^{ c}} \DD_{y^{k}_{v}/\sqrt{\la}}\)
\( \prod_{ i\in B_{s}}\mathcal{D}_{ x^{ s( i)}_{ c(
i)}/\sqrt{\la}}\)\label{4h.7}\\ &&
\prod_{ i\in B^{ c}_{ s}}G_{\la}\(z_{s( i) }+\sum_{ j=2}^{ c( i)}
x^{  s( i)}_{ j}/\sqrt{\la} -(z_{s( i-1) }+\sum_{ j=2}^{ c( i-1)} x^{ s(
i-1)}_{ j}/\sqrt{\la} )\)\Bigg|\nn\\ &&
\leq c     \la^{ -2}\prod_{k\in F^{ c}}\,|y^{k}_{v} |^{\nu} \prod_{ 
i\in B_{ s}}\,|x^{ s( i)}_{ c( i)} |^{\nu}.\nn
\end{eqnarray} Our lemma then follows using (\ref{4h.5a}) which
implies that
$|F^{ c}|+| B_{s} |\geq 2$. \qed

\section{Approximating intersection local
times}\label{sec-approxilt} 

The goal of this section is to prove the following lemma.

\begin{lemma}\label{lem-rwapprox} We can find a Brownian
motion such for each $j\geq 1$ there exists
$\beta>0$ such that
\begin{equation}
\|\la\Ga_{j,\la}( \ze_{  \la})-\ga_{j}( \ze,\om_{\la^{ -1}})\|_{
2}=O(\la^{
\beta} )\label{32.1}
\end{equation}
\end{lemma}

{\bf Proof of Lemma \ref{lem-rwapprox}:} Let $f( x)$ be a smooth
function  on
$R^2$, supported in the unit disc and with $\int f( x)\,dx=1$. We
set
$f_{
\ep}( x)={  1\over \ep^{ 2}}f( x/\ep)$. On the one hand it is easy
to  see that if we set
$\wt u^{ 1}( f_{ \tau})=\int u^{ 1}( x)f_{ \tau}( x)\,dx$ and
\begin{eqnarray}
\wt{\ga}_{k}(\ze,f_{ \tau})&=&\int \ol\ga_{k}(\ze,x)\,\prod_{ i=2}^{
k}f_{
\tau}( x_{ i})\,dx_{ 2} \cdots dx_k, \nn\\
\wt{\al}_{j}(\ze,f_{ \tau})&=&\int \ol\al_{j}(\ze,x)\,\prod_{  i=2}^{
j}f_{
\tau}( x_{ i})\,dx_{ 2} \cdots dx_k,\nn
\end{eqnarray}
     we will have
\begin{equation}
\wt{\ga}_{k}(\ze,f_{ \tau})=\sum_{j=1}^{ k}{k-1 \choose j-1}
\(-\wt u^{ 1}( f_{ \tau})\)^{k- j}\wt{\al}_{j}(\ze,f_{ \tau})\label{32.2}
\end{equation} and
\begin{equation}
\wt{\al}_{j}(t,f_{ \tau})=\int_{ 0\leq t_{ 1}\leq\cdots\leq t_{ j}<t}
\prod_{ i=2}^{ j} f_{ \tau}(W_{ t_{ i}}-W_{ t_{ i-1}})dt_{ 1}\cdots
dt_{ j}.\label{32.3}
\end{equation} On the other hand it follows from (\ref{7.6}) and
Jensen's inequality that
\begin{equation}
\|\wt{\ga}_{k}(\ze,f_{ \tau})-\ga_{k}(\ze) \|_{ 2}\leq C\tau^{
\bar{\de}}.\label{32.4}
\end{equation}

If we set
$\wt G_{ \la}( f_{ \tau})=\int G_{ \la}( x/\sqrt{\la})f_{ \tau}( x)\,dx$,
\[
\wt{\Ga}_{k,\la}(\ze_{ \la},f_{ \tau})=\int \ol \Ga_{k,\la}(\ze_{
\la},x)\,\prod_{ i=2}^{ k}f_{
\tau}( x_{ i})\,dx_{ 2} \cdots dx_k\] and
\[\wt{I}_{j}(\ze_{ \la},f_{ \tau})=\int \ol I_{j}(\ze_{
\la},x/\sqrt{\la})\,\prod_{ i=2}^{ j}f_{
\tau}( x_{ i})\,dx_{ 2} \cdots dx_k,\]
     we similarly have
\begin{equation}
\wt{\Ga}_{k,\la}(\ze_{ \la},f_{ \tau})=\sum_{j=1}^{ k}{k-1 \choose
j-1}
\(-\wt G_{ \la}( f_{ \tau})\)^{k- j}\wt{I}_{j}(\ze_{ \la},f_{
\tau})\label{32.5}
\end{equation} and
\begin{eqnarray}\qquad \wt{I}_{j}(\ze_{\la},f_{ \tau}) &=& \int_{
0\leq t_{ 1}\leq\cdots\leq t_{  j}<\ze_{\la}} \prod_{ i=2}^{ j} f_{
\tau}(\sqrt{\la}(S_{ [t_{ i}]}-S_{ [t_{ i-1}]}))dt_{  1}\cdots dt_{
j}\label{32.6}\\ &=& \int_{ 0\leq t_{ 1}\leq\cdots\leq t_{ j}<\ze}
\prod_{ i=2}^{ j} f_{
\tau}(\sqrt{\la}(S_{ [t_{ i}/\la]}-S_{ [t_{  i-1}/\la]}))dt_{ 1}\cdots
dt_{ j}.\nonumber
\end{eqnarray} 
 It then follows from (\ref{7.7}) that
\begin{equation}
\| \la \wt{\Ga}_{k,\la}( \ze_{\la},f_{ \tau})-\la\Ga_{k,\la}(
\ze_{\la})\|_{ 2}\leq C\tau^{
\bar{\de}}.\label{32.7}
\end{equation}

To complete the proof of Lemma \ref{lem-rwapprox} it only
remains to show that with $\tau=\la^{ \rho}$ for $\rho>0$ small
\begin{equation}
\|\la \wt{\Ga}_{k,\la}( \ze_{\la},f_{ \tau})-\wt{\ga}_{k}(\ze,f_{
\tau},\om_{
\la^{ -1}}) \|_{ 2}
\leq c\la^{ \ze}\label{4.8}
\end{equation} for some $c<\ff$ and $\ze>0$.
By (\ref{32.2})-(\ref{32.6}) it suffices to show that for some 
$\de'>0$ and all sufficiently small $\tau,\la$
\begin{equation}
\wt u^{ 1}( f_{ \tau})
=O( \log ( 1/|\tau|)),\hspace{ .3in}
|\wt G_{ \la}( f_{ \tau})-\wt u^{ 1}( f_{ \tau})|
\leq c\tau^{-3}\la^{ \de'}\label{4.8m}
\end{equation}
and
\begin{eqnarray}
\|\wt{\al}_{k}(\ze,f_{\tau},\om_{\la^{ -1}}) \|_{ 2}
& \leq & c\tau^{-2(k-1)},\nn\\
\|\la \wt{I}_{k,\la}( \ze_{\la},f_{ \tau})-\wt{\al}_{k}(\ze,f_{
\tau},\om_{
\la^{ -1}}) \|_{ 2}
& \leq  & c\tau^{-2k+1}\la^{ \de'}.\label{4.8n}
\end{eqnarray}

The first part of (\ref{4.8m}) follows from the fact that
 $u^{ 1}( x)=O( \log ( 1/|x|))$, see 
\cite[(2.b)]{LG90}.
To prove the second part of (\ref{4.8m}), we note that $\sup_{  x}
|\nabla f_{\tau}(x)|\leq c\tau^{- 3}$, so
\begin{eqnarray} &&
|\wt G_{ \la}( f_{ \tau})-\wt u^{ 1}( f_{ \tau})|\label{4.9k}\\ && 
=\left|\int_{ 0}^{ \ff}e^{ -t}E\(f_{ \tau}(  \sqrt{\la}  S_{
[t/\la]})-f_{ \tau}(\sqrt{\la} W_{t/\la})\) dt\right|\nonumber\\ && 
\leq c\tau^{- 3}\int_{ 0}^{ \ff}e^{ -t}\Bigg\|   \sqrt{\la}  (S_{
[t/\la]}-W_{t/\la})\Bigg\|_{ 1} dt.\nonumber
\end{eqnarray}
The second part of (\ref{4.8m}) then follows from
 the last inequality in Section
\ref{sec-strongapp}.

The first part of (\ref{4.8n}) follows from the fact that
 $\sup_{  x}| f_{\tau}(x)|\leq c\tau^{- 2}$, so
that
\begin{equation}
\|\wt{\al}_{k}(\ze,f_{\tau},\om_{\la^{ -1}}) \|_{ 2}^{ 2}
\leq c\tau^{-2(k-1)}\int_{ 0}^{ \ff}e^{ -t}t^{ n} dt.\label{4.9p}
\end{equation}

To prove the second part of (\ref{4.8n}), we use the above bounds
on 
$\sup_{  x} |\nabla f_{\tau}(x)|$ and  $\sup_{  x}|
f_{\tau}(x)|$ to see that
\begin{eqnarray} &&\qquad
\|\la \wt{I}_{k,\la}( \ze_{\la},f_{ \tau})-\wt{\al}_{k}(\ze,f_{
\tau},\om_{
\la^{ -1}}) \|_{ 2}^{ 2}\label{4.9q}\\  && 
\leq c\tau^{-2k+1}\sum_{ j=1}^{ k}\int_{ 0}^{ \ff}e^{ -t}
\(\int_{ 0\leq t_{ 1}\leq\cdots\leq t_{ k}<t}
\Bigg\|   \sqrt{\la}  (S_{
[t_{ j}/\la]}-W_{t_{ j}/\la})\Bigg\|^{ 2}_{ 2}dt_{ 1}\cdots dt_{ k}\)
dt.\nonumber
\end{eqnarray}
The second part of (\ref{4.8n}) then follows from
 the last inequality in Section
\ref{sec-strongapp}.

\section{Renormalized Brownian intersection local
times}\label{sec-brilt}

Recall the definition of $\gamma_k(t)$ given in (\ref{2.4a}). Note
  from \cite[(2.b)]{LG90} that for some fixed constant c
\begin{equation} u_{ \ep}=\int_{ 0}^{ \ff} e^{ -t}p_{t+ \ep}(0)
\,dt={1
\over 2\pi}\log ( 1/\ep) +c+O( \ep).\label{2.11}
\end{equation} In \cite{RosenJC} we show that the limit in
(\ref{2.4a}) exists a.s.  and in all $L^{ p}$ spaces, and that
$\ga_{k}(t)$ is continuous  in $t$. The rest of this section  is
basically contained in
\cite{LG90} but we point out that \cite{RosenJC} came after
\cite{LG90} and resulted in some simplification.

For any given function $h:(0,\infty)\to R$ we set
$\wh \ga_{ 1}(t,h)=t$ and  for $k\geq 2$
\begin{equation}\quad
\wh \ga_{k}(t,h)=\lim_{ \ep\rar 0}\sum_{ l=1}^{ k}{k-1\choose l-1} (
-h_{
\ep})^{ k-l}
\alpha_{l,\ep}(t),\label{2.10h}
\end{equation}  where we write $h_\ep$ for $h(\ep)$. In
particular,
$\ga_{ 1}(t)=\wh \ga_{ 1}(t,u)$.
     Let $\mathcal{H}$ denote the set of functions $h$ such that
     $\lim_{ \ep\rar 0}(h_{ \ep}-u_{ \ep})$ exists and is finite. In 
the next lemma we will see that the limit in  (\ref{2.10h}) exists
for all
$h\in
\mathcal{H}$.

\begin{lemma}[Renormalization Lemma]\label{lem-renorm} Let
$h\in
\mathcal{H}$. Then $\wh \ga_{k}(t,h)$ exists for all $k\geq 1$ and if
$\bar{h}\in
\mathcal{H}$ with
$\lim_{ \ep\rar 0}(h_{ \ep}-\bar{h}_{ \ep})=b$ then for any
$k\geq 1$
\begin{equation}
\wh \ga_{k}(t,h)=\sum_{ m=1}^{ k}{k-1\choose m-1}(-b)^{k-m}\wh \ga_{
m}(t,\bar{h}).\label{2.11h}
\end{equation}
\end{lemma} {\bf Proof of Lemma \ref{lem-renorm}:}  Setting $b_{
\ep}=h_{ \ep}-\bar{h}_{ \ep}$ we have
\begin{eqnarray} &&
\sum_{ l=1}^{ k}{k-1\choose l-1} ( -h_{ \ep})^{ k-l}\al_{ l, \ep}(t
)\label{2.14h}\\ &&
     =\sum_{ l=1}^{ k}{k-1\choose l-1} (-\bar{h}_{ \ep}-b_{ \ep})^{
k-l}\al_{ l,
\ep}(t )\nonumber\\ &&
     =\sum_{ l=1}^{ k}{k-1\choose l-1}\sum_{ j=0}^{
k-l}{k-l\choose j} (-b_{
\ep})^{ j}(-\bar{h}_{ \ep})^{ (k-j)-l}\al_{ l, \ep}(t )\nonumber
\end{eqnarray} Using
\[{k-1\choose l-1}{k-l\choose j}={k-1\choose j}{k-j-1\choose
l-1}\]  the last line in (\ref{2.14h}) becomes
\be
\sum_{ j=0}^{ k-1}{k-1\choose j}(-b_{ \ep})^{ j}\sum_{ l=1}^{
k-j}{k-j-1\choose l-1} (-\bar{h}_{ \ep})^{ (k-j)-l}\al_{ l, \ep}(t
).\label{2.15h}
\ee Taking $\bar{h}_{ \ep}=u_{ \ep}$ then shows the existence of
$\ga_{k}(t,h)$. Returning to general  $\bar{h}\in \mathcal{H}$
and now taking the $\ep\rar 0$  limit we obtain
     \begin{eqnarray} &&
\wh \ga_{ k}(t,h)=\sum_{ j=0}^{ k-1}{k-1\choose j}(-b)^{ j}\wh \ga_{
k-j}(t,\bar{h} )\label{2.16h}\\ &&
\hspace{ .6in}=\sum_{ m=1}^{ k}{k-1\choose m-1}(-b)^{k-m}\wh \ga_{
m}(t,\bar{h})\nonumber
\end{eqnarray} where the last line follows from the substitution
$m=k-j$.
\qed

  Let $h\in \mathcal{H}$. We shall sometimes write
$\wh \ga_{k}(t,h,\om)$ for
$\wh \ga_{k}(t,h)$ to emphasize its dependence on the path $\om$.
We want to discuss how renormalized intersection local time
changes with a time rescaling. Let
$\om_{ r}( s)=r^{ -1/2}\om( rs)$. 
Then $\wh \ga_k(t,h,\om_r)$ is the same as $\wh \ga_k(t,h)$ defined
in terms of the Brownian motion $W^{(r)}_t=W_{rt}/\sqrt r$.

\begin{lemma}[Rescaling Lemma]\label{lem-rescale}  Let $h\in
\mathcal{H}$.  Then for any
$k\geq 1$
\begin{equation} 
\wh \ga_{k}(t,h,\om_{ r})=r^{ -1}\sum_{ m=1}^{ k}{k-1\choose
m-1}\({1
\over 2\pi}\log ( 1/r)\)^{k-m}\wh \ga_{ m}(rt,h,\om).\label{2.11s}
\end{equation}
\end{lemma}  {\bf Proof of Lemma \ref{lem-rescale}:} After
replacing
$\om$ by $\om_{ r}$  the integral on the right hand side of
(\ref{2.10h}) is replaced by
\begin{eqnarray} &&
\int_{ 0\leq t_{ 1}\leq\cdots\leq t_{ l}<t} \prod_{ i=2}^{ l} p_{
\ep}\({W_{ rt_{ i}}-W_{ rt_{ i-1}} \over \sqrt{r}}\)dt_{ 1}\cdots 
dt_{ l}\label{2.12}\\ &&=r^{-l}\int_{ 0\leq t_{ 1}\leq\cdots\leq
t_{ l}<rt}
\prod_{ i=2}^{ l} p_{ \ep}\({W_{ t_{ i}}-W_{ t_{ i-1}} \over
\sqrt{r}}\)dt_{ 1}\cdots dt_{ l}
\nonumber\\ &&=r^{-1}\int_{ 0\leq t_{ 1}\leq\cdots\leq t_{ l}<rt}
\prod_{ i=2}^{ l} p_{ r\ep}(W_{ t_{ i}}-W_{ t_{ i-1}} )dt_{ 1}\cdots
dt_{ l}.
\nonumber
\end{eqnarray} Abbreviating this last integral as $\al_{ l,
r\ep}(rt,\om )$ we have
\begin{equation}
\wh \ga_{k}(t,h,\om_{ r})=r^{-1}\lim_{ \ep\rar 0}\sum_{ l=1}^{
k}{k-1\choose l-1} ( -h_{ \ep})^{ k-l}\al_{ l, r\ep}(rt,\om
).\label{2.13}
\end{equation} Since $h\in \mathcal{H}$ it is easily seen that
$\lim_{ \ep\rar 0}(h_{ \ep}-h_{ r\ep})=-{1 \over 2\pi}\log (
1/r)$ and our lemma then follows from Lemma \ref{lem-renorm}.
\qed

\section{Range and Brownian intersection local
times}\label{sec-rabilt} In this section we prove the following
theorem. 
\bt \label{theo-expon} For each $k\geq 1$
\begin{equation} g^{k}_{ \la}
\(\la |\mathcal{R}( \ze_{  \la })|-\sum_{ j=1}^{k}( -1)^{ j-1}g^{
-j}_{
\la }
\ga_{j}( \ze,\om_{\la^{ -1}})\)\rar 0\hspace{ .2in}a.s.\label{2.9}
\end{equation} as $\la\to 0$.
\et

{\bf Proof of Theorem \ref{theo-expon}:} Using (\ref{32.1})
together with Lemma \ref{lem-rwasip} and its  proof, we see that
for some
$M_{k}<\ff$
\begin{equation}
     \Bigg\| g^{4k+1 }_{ \la}
\(\la|\mathcal{R}( \ze_{  \la})|-\sum_{ j=1}^{4k}( -1)^{ j-1}g^{
-j}_{
\la}
\ga_{j}( \ze,\om_{\la^{ -1}})\)\Bigg\|_2^{ 2}\leq M_{
k}\label{2.1}
\end{equation} for all $\la>0$ sufficiently small.

We now follow Le Gall \cite{LG90}. With $\la_{ n}=e^{ -n^{
1/2k}}$ we have that for any $\ep>0$
\begin{eqnarray} &&
\sum_{ n=1}^{ \ff}P\lc g^{k }_{ \la_{ n}}
\(\la_{ n}|\mathcal{R}( \ze_{  \la_{ n}})|-\sum_{ j=1}^{4k}( -1)^{
j-1}g^{ -j}_{
\la_{ n}}
\ga_{j}( \ze,\om_{\la_{ n}^{ -1}})\)\geq g^{-1}_{ \la_{
n}}\rc\label{2.2}\\ &&
\leq \sum_{ n=1}^{ \ff}P\lc g^{4k+1 }_{ \la_{ n}}
\(\la_{ n}|\mathcal{R}( \ze_{  \la_{ n}})|-\sum_{ j=1}^{4k}( -1)^{
j-1}g^{ -j}_{
\la_{ n}}
\ga_{j}( \ze,\om_{\la_{ n}^{ -1}})\)\geq g^{3k}_{ \la_{ n}}\rc
\nonumber\\ &&
\leq M_{ k}\sum_{ n=1}^{ \ff}  g^{-6k }_{ \la_{ n}}<\ff.
\nonumber
\end{eqnarray}
     Then by Borel-Cantelli
\begin{equation} g^{k }_{ \la_{ n}}
\(\la_{ n}|\mathcal{R}( \ze_{  \la_{ n}})|-\sum_{ j=1}^{4k}( -1)^{
j-1}g^{ -j}_{
\la_{ n}}
\ga_{j}( \ze,\om_{\la_{ n}^{ -1}})\)\rar 0\hspace{
.2in}a.s.\label{2.3aa}
\end{equation}  Since  for each $m\geq 1$ we have that
$\ga_{j}(\ze,\om_{\la_{ n}^{ -1}})$ is bounded in $L^m$
uniformly in
$n$, then by Chebyshev's inequality with $m$ sufficiently large
$P(\ga_j(\ze, \om_{\la_n^{-1}})> g_{\la_n})$  will be summable.
So we may drop the terms for $j>k$ and we then have
\begin{equation} g^{k }_{ \la_{ n}}
\(\la_{ n}|\mathcal{R}( \ze_{  \la_{ n}})|-\sum_{ j=1}^{k}( -1)^{ 
j-1}g^{ -j}_{
\la_{ n}}
\ga_{j}( \ze,\om_{\la_{ n}^{ -1}})\)\rar 0\hspace{
.2in}a.s.\label{2.3}
\end{equation}

Before continuing the proof of Theorem \ref{theo-expon} we first
prove the following lemma.
\begin{lemma}\label{lem-hct} For  any $k\geq 1$
\begin{equation}
\lim_{ n\rar 0}\sup_{\la_{ n+1}\leq  \la\leq \la_{ n}}
     |\ga_{k}( \ze,\om_{\la^{ -1}})-\ga_{k}(
\ze,\om_{\la_{ n}^{ -1}})|=0. ~~~~a.s.\label{2.01}
\end{equation}
\end{lemma}

{\bf Proof of Lemma \ref{lem-hct}:} By (\ref{2.11s})
     for any $k\geq 1$
\begin{equation}
\ga_{k}(\ze,\om_{\la^{ -1}}) ={\la \over \la_{ n}}\sum_{ m=1}^{
k}{k-1\choose m-1}\({1 \over 2\pi}\log \Big( {\la \over \la_{
n}}\Big)\)^{k-m}\ga_{ m}\Big( {\la_{ n} \over
\la}\ze,\om_{\la_{ n}^{ -1}}\Big).\label{2.02}
\end{equation} Hence for any $p\geq 1$
\begin{eqnarray} &&
\|   \sup_{\la_{ n+1}\leq  \la\leq \la_{ n}}
     |\ga_{k}( \ze,\om_{\la^{ -1}})-\ga_{k}(
\ze,\om_{\la_{ n}^{ -1}})|\, \|_{ p}\label{2.03}\\ &&
\leq \|   \sup_{\la_{ n+1}\leq  \la\leq \la_{ n}}
     |{\la \over \la_{ n}}\ga_{ k}( {\la_{ n} \over
\la}\ze,\om_{\la_{ n}^{ -1}})-\ga_{k}(
\ze,\om_{\la_{ n}^{ -1}})|\, \|_{ p}
     \nonumber\\ &&  + c\sum_{ m=1}^{ k-1}\sup_{\la_{ n+1}\leq
\la\leq
\la_{ n}}
\({1 \over 2\pi}\log ( {\la \over \la_{ n}})\)^{k-m}
\|\sup_{\la_{ n+1}\leq  \la\leq \la_{ n}}\ga_{ m}( {\la_{ n} \over
\la}\ze,\om_{\la_{ n}^{ -1}})\|_{ p} \nn\\ && = \|   \sup_{\la_{
n+1}\leq
\la\leq \la_{ n}}
     |{\la \over \la_{ n}}\ga_{ k}( {\la_{ n} \over \la}\ze)-\ga_{k}(
\ze)| \|_{ p}
     \nonumber\\ &&  + c\sum_{ m=1}^{ k-1}\sup_{\la_{ n+1}\leq
\la\leq
\la_{ n}}
\({1 \over 2\pi}\log ( {\la \over \la_{ n}})\)^{k-m}
\|\sup_{\la_{ n+1}\leq  \la\leq \la_{ n}}\ga_{ m}( {\la_{ n} \over
\la}\ze)\|_{ p} \nn
\end{eqnarray} It follows from (9.11) of \cite{BK} that for any
$k\geq 1$ we can  find $\bb>0$ such that
\begin{equation}
\| \sup_{\stackrel{|t-s|\leq \de}{s,t\leq 1}}|\ga_{ k}(s)-\ga_{k}(t)|
\,\|_{p}\leq c\de^{\bb}.\label{2.04}
\end{equation} Actually, this is proved  for a renormalized
intersection local time
${\xi}_{ k}(t)$ where ${\xi}_{ k}(t)=\lim_{ x\rar
0}{\xi}_{ k}(t,x)$ and
${\xi}_{ k}(t,x)$ differs from $\ol\ga_{ k}(t,x)$ defined in 
(\ref{7.5}) in that
$u^{ 1}(  x)$ is replaced by
$\pi^{-1}\log (1 / |x|)$. Since
$u^{ 1}( x)- \pi^{-1}\log ({1 / |x|})=c+O(|x|^{2}\log |x|)$, see
\cite[(2.b)]{LG90}, we obtain  (\ref{2.04}). Using (\ref{2.11s})
again we find that
\begin{equation}\hspace{ .5in}
\|   \sup_{\la_{ n+1}\leq  \la\leq \la_{ n}} |\ga_{ k}( {\la_{ n}
\over
\la}t)-\ga_{k}( t)|\, \|_{ p}\leq ct( \log t)^{ k}\Bigg| {\la_{ n}
\over\la_{ n+1}}-1\Bigg|^{\bb}
\leq ct( \log t)^{ k}n^{-\bb'}
\label{2.05}
\end{equation} where we have used
\begin{equation}
\log {\la_{ n} \over \la_{ n+1}}=O(n^{-1+ 1/k} ).\label{2.06a}
\end{equation} Hence
\begin{equation}
\|   \sup_{\la_{ n+1}\leq  \la\leq \la_{ n}} |\ga_{ k}\Big( {\la_{ n}
\over
\la}\ze\Big)-\ga_{k}(\ze)| \, \|_{ p}
\leq cn^{-\bb''}.
\label{2.06}
\end{equation} Using (\ref{2.03}) and (\ref{2.06}) now shows that
\begin{equation}
\|   \sup_{\la_{ n+1}\leq  \la\leq \la_{ n}}
     |\ga_{k}( \ze,\om_{\la^{ -1}})-\ga_{k}(
\ze,\om_{\la_{ n}^{ -1}})|  \,\|_{ p}  \leq cn^{-\bb''}\label{2.07}
\end{equation} and our lemma then follows using Holder's
inequality for sufficiently large $p$ and the Borel-Cantelli Lemma.
\qed

Continuing the proof of Theorem \ref{theo-expon}, by our choice
of
$ \la_{ n}$
\begin{equation}
\lim_{ n\rar 0}g^{k}_{ \la_{ n+1}}-g^{k}_{ \la_{ n}}=0.\label{2.6}
\end{equation} Together with (\ref{2.01}) we have that a.s.
\begin{equation}\qquad
\lim_{ n\rar 0}\sup_{\la_{ n+1}\leq  \la\leq \la_{ n}} |\sum_{
j=1}^{k}( -1)^{ j-1}g^{k -j}_{ \la}
\ga_{j}( \ze,\om_{\la^{ -1}})-\sum_{ j=1}^{k}( -1)^{ j-1}g^{k -j}_{
\la_{ n}}
\ga_{j}( \ze,\om_{\la_{ n}^{ -1}})|=0.\label{2.7}
\end{equation}

Using the fact that $|\mathcal{R}( \ze_{  \la})|$ and $g_{ \la}$ are
monotone decreasing we have that
\begin{eqnarray} &&
\sup_{\la_{ n+1}\leq  \la\leq \la_{ n}}
\Big|\la g^{k}_{ \la}|\mathcal{R}( \ze_{  \la})|-
\la_{ n}g^{k}_{ \la_{ n}}|\mathcal{R}( \ze_{  \la_{
n}})|\Big|\label{2.8}\\ &&
\leq \Big|\la_{ n} g^{k}_{ \la_{ n+1}}|\mathcal{R}( \ze_{  \la_{
n+1}})|-
\la_{ n+1} g^{k}_{ \la_{ n}}|\mathcal{R}( \ze_{  \la_{ n}})|\Big|
\nonumber\\ &&
\leq \Big|\la_{ n} g^{k}_{ \la_{ n+1}}|\mathcal{R}( \ze_{  \la_{
n+1}})|-
\la_{ n+1} g^{k}_{ \la_{ n+1}}|\mathcal{R}( \ze_{  \la_{
n+1}})|\Big|
\nonumber\\ && +\Big|\la_{ n+1} g^{k}_{ \la_{
n+1}}|\mathcal{R}(
\ze_{  \la_{ n+1}})|-
\la_{ n} g^{k}_{ \la_{ n}}|\mathcal{R}( \ze_{  \la_{ n}})|\Big|
\nonumber\\ && + \Big|\la_{ n} g^{k}_{ \la_{ n}}|\mathcal{R}(
\ze_{
\la_{ n}})|-
\la_{ n+1} g^{k}_{ \la_{ n}}|\mathcal{R}( \ze_{  \la_{ n}})|\Big|
\nonumber\\  &&
\leq 2|\la_{ n}-\la_{ n+1}|g^{k}_{ \la_{ n+1}}|\mathcal{R}( \ze_{
\la_{ n+1}})|\nonumber\\ && + \Big|\la_{ n+1} g^{k}_{ \la_{
n+1}}|\mathcal{R}(
\ze_{  \la_{ n+1}})|-
\la_{ n}g^{k}_{ \la_{ n}}|\mathcal{R}( \ze_{  \la_{ n}})|\Big|\rar
0\hspace{ .2in}a.s.\nn
\end{eqnarray} Here the first term on the right hand side of 
(\ref{2.8}) goes to
$0$  using the fact that
\[|\la_{ n}-\la_{ n+1}|=|1-e^{n^{1/2k}-(n+1)^{ 1/2k} } |\la_{
n}\leq  n^{ -1}\la_{ n}
\leq 2n^{ -1}\la_{ n+1}\] and the second term on the right hand
side of (\ref{2.8}) goes to $0$  using (\ref{2.7}) and (\ref{2.3}).
Combining (\ref{2.3}), (\ref{2.7}) and (\ref{2.8}) we have
(\ref{2.9}).
\qed

\section{Non-random times}\label{sec-nonrandom}

In this section we complete the proof of Theorem
\ref{theo-fixasip}.  Recall that $\zeta_\la=n$ if
$n-1<\frac{1}{\la}\zeta\leq n$.  So $\zeta_\la=
\lceil \frac{1}{\la}\zeta\rceil$ where $\lceil x\rceil$ denotes the
smalllest integer $m\geq x$. Hence (\ref{2.9}) can be written as
\begin{equation} \label{2.9bx}
g_{\la}^{2k}\(\la|\mathcal{R}(\lceil
\ze/\la\rceil )| - \sum_{j=1}^{2k} (-1)^{j-1} g_{\la}^{-j} \ga_j(\ze,
\om_{\la^{-1}})\)
\to 0, \qquad a.s.
\end{equation}

If $(\Om, P)$ denotes our probability space for
$\{S_{ n}\,;\,n\geq 1\}$ and
$\{W_{ t}\,;\,t\geq 0\}$, then the almost sure convergence in
(\ref{2.9bx}) is with respect to the measure $e^{ -t}\,dt\times P$
on
$R^{ 1}_{ +}\times \Om$, where $\ze(t,\om)=t$. Hence by Fubini's
theorem we have that for almost every $t>0$
\begin{equation} \label{2.9x} 
g_{\la}^{2k}\(\la|\mathcal{R}(\lceil
t/\la\rceil )| - \sum_{j=1}^{2k} (-1)^{j-1} g_{\la}^{-j} \ga_j(t,
\om_{\la^{-1}})\)
\to 0, \qquad a.s.
\end{equation}

Fix a $t_0$ for which (\ref{2.9x}) holds and let $\la$ run through
the sequence
$t_0/n$. Then (\ref{1.5}) and (\ref{2.9x}) tell us that
\begin{equation} (\log n)^{k}\(\frac{t_0}{n}
|\mathcal{R}(n)|+\sum_{j=1}^{k}  (-g_{t_0/n})^{-j}
\gamma_j(t_0,\om_{n/t_0})\)\to 0 \hspace{ .2in}
a.s.\label{b.8}
\end{equation} Using (\ref{2.11s}) and writing $b_{ r}={1 \over
2\pi}\log ( 1/r)$ we have that
\begin{equation} \qquad (\log n)^{k}\(\frac{t_0}{n}
|\mathcal{R}(n)|+t_0\sum_{j=1}^{k} (-g_{t_0/n})^{-j}
\sum_{ m=1}^{ j}{j-1\choose m-1}b_{1/t_{ 0}}^{j-m}\ga_{
m}(1,\om_{ n})\)\to 0
\hspace{ .2in} a.s.\label{b.9}
\end{equation}

Then
\begin{eqnarray} &&
\sum_{j=1}^{k} (-g_{t_0/n})^{-j}
\sum_{ m=1}^{ j}{j-1\choose m-1}b_{1/t_{ 0}}^{j-m}\ga_{
m}(1,\om_{ n})\label{b.11}\\ &&  =\sum_{m=1}^{k}
\(\sum_{ j=m}^{ k}{j-1\choose m-1}\({-b_{1/t_{ 0}} \over
g_{t_0/n}}\)^{j-m}\) (-g_{t_0/n})^{-m}\ga_{ m}(1,\om_{
n}).\nonumber
\end{eqnarray}

Now,
\begin{equation}\qquad
\sum_{ j=m}^{ k}{j-1\choose m-1}x^{j-m}=\sum_{ i=0}^{
k-m}{i+m-1\choose m-1}x^{i}=\({ 1\over 1-x}\)^{ m}+O(x^{k-m+1}
).\label{b.12}
\end{equation}

By (\ref{2.04}) with $\de=1$ we have that
$\sup_{t\leq 1} |\gamma_j(t,\om)|$ is in $L^p$ for each $p$ and
each
$j\geq 1$. If we set $V_{j,\ell}=\sup_{t\leq 1} |\gamma_j(t,
\om_{2^\ell})|$, we then have, taking $p$ large enough, that
\[
\sum_{\ell=1}^\infty P(V_{j,\ell}> \eta \log(2^\ell))
\leq \sum_{\ell=1}^\infty \frac{E V_{j,\ell}^p}{(\eta \log
2^\ell)^p}
\]
  is summable for each $\eta$. Hence by Borel-Cantelli
$V_{j,\ell}/\log(2^\ell)\to 0$ a.s.  for each $j\geq 1$. Since by
Lemma
\ref{lem-rescale} we have for $2^\ell\leq r< 2^{\ell+1}$
   that $\ga_k(1,\om_r)$ is bounded by a linear combination   of
the
$V_{j,\ell}$,
$1\leq j\leq k$, with coefficients that are bounded independently
of
$r$, we conclude
\[
\ga_j(1, \om_n)/\log n \to 0, \qquad a.s.
\] Thus we can replace (\ref{b.11}) up to errors which are
$O\(\log n
\)^{ -k-1}$ by
\begin{eqnarray} &&
\sum_{m=1}^{k}
\({ -1\over g_{t_0/n}+b_{1/t_{ 0}}}\)^{ m}
\ga_{ m}(1,\om_{ n})=\sum_{m=1}^{k} (-g_{1/n})^{-m}
\ga_{ m}(1,\om_{ n})\label{b.13}
\end{eqnarray} since by (\ref{1.5}) we have that
$g_{t_0/n}+b_{1/t_{ 0}}=g_{1/n}+O(n^{-\de} )$.

     Thus we obtain
\begin{equation} (\log n)^{k}\(\frac{1}{n}
|\mathcal{R}(n)|+\sum_{j=1}^{k} (-g_{1/n})^{-j} \ga_{ j}(1,\om_{
n})\)\to 0
\hspace{ .2in} a.s.\label{b.14}
\end{equation}

This, together with (\ref{45.1}), gives Theorem \ref{theo-fixasip}.
\qed

\appendix
\section{Estimates for random walks}\label{sec-estimates}

In this appendix we will obtain some estimates for strongly
aperiodic planar random walks
$S_n=\sum_{i=1}^n X_i$, where the $X_i$ are symmetric, have the
identity as covariance matrix, and have $2+\de$ moments for
some
$\de>0$.

Let
\[G_{ \la}(x){:=}\sum_{n=0}^{\ff} e^{ -\la n} q_{n}(x).\] If
\[
\phi(p)=E(e^{ip\cdot X_1})
\] denotes the characteristic function of $X_1$, we have
\be G_{ \la}(x)={1\over (2\pi)^2}\int_{[-\pi,\pi]^2} {e^{ ip\cdot
x}\over 1-e^{ -\la }\phi(p)}\,dp.\label{45.0}
\ee

\begin{lemma}\label{lem-green} Let $S_{n}$ be as above.  Then
\be
     G_{ \la}(0)={1 \over 2\pi}\log ( 1/\la)+c_{ X} +O( \la^{ \de
}\log (1/\la))\label{45.1}
\ee where
\begin{equation} c_X={1 \over 2\pi}\log (\pi^{ 2}/2 )+ {1 \over
(2\pi)^{ 2}}\int_{[-\pi,\pi]^2} {\phi(p)-1+|p|^{ 2}/2\over
(1-\phi(p))|p|^{ 2}/2}\,dp\label{1.5cy}
\end{equation} is a finite constant.
\end{lemma}

{\bf Proof of Lemma \ref{lem-green}:}  We have
\be G_{ \la}(0)={1\over (2\pi)^2}\int_{[-\pi,\pi]^2} {1\over
1-e^{ -\la }\phi(p)}\,dp.\label{45.3}
\ee We intend to compare this with
\[{1\over (2\pi)^2}\int_{[-\pi,\pi]^2} {1\over
\la+|p|^{ 2}/2}\,dp\] whose asymptotics are easier to compute.
Indeed,
\bea &&
\int_{[-\pi,\pi]^2} {1\over
\la+|p|^{ 2}/2}\,dp\label{45.13}\\ &&=\int_{D( 0,\pi)} {1\over
\la+|p|^{ 2}/2}\,dp+\int_{[-\pi,\pi]^2-D( 0,\pi)} {1\over
\la+|p|^{ 2}/2}\,dp\nn
\eea where $D( 0,\pi)$ is the disc centered at the origin of radius
$\pi$. It is clear that
\begin{equation}
\int_{[-\pi,\pi]^2-D( 0,\pi)} {1\over
\la+|p|^{ 2}/2}\,dp=\int_{[-\pi,\pi]^2-D( 0,\pi)} {1\over |p|^{
2}/2}\,dp+O(
\la).\label{45.14}
\end{equation} On the other hand, using polar coordinates
\begin{equation}
\int_{D( 0,\pi)} {1\over
\la+|p|^{ 2}/2}\,dp=2\pi\(\log (\la+\pi^{ 2}/2 )-\log
(\la)\).\label{45.15}
\end{equation} Thus
\begin{equation} {1\over (2\pi)^2}\int_{[-\pi,\pi]^2} {1\over
\la+|p|^{ 2}/2}\,dp={1 \over 2\pi}\log ( 1/\la) +{1 \over 2\pi}\log
(\pi^{ 2}/2 )+O(\la).\label{45.15r}
\end{equation}

We then note that
\begin{eqnarray} &&
\int_{[-\pi,\pi]^2} {1\over 1-e^{ -\la 
}\phi(p)}\,dp-\int_{[-\pi,\pi]^2} {1\over
\la+|p|^{ 2}/2}\,dp\label{45.4}\\ && =\int_{[-\pi,\pi]^2} {(
\la+|p|^{ 2}/2)-(1-e^{ -\la }\phi(p))\over (1-e^{ -\la }\phi(p))(
\la+|p|^{ 2}/2)}\,dp\nonumber\\ && =\int_{[-\pi,\pi]^2}
{\phi(p)-1+|p|^{ 2}/2\over (1-e^{ -\la }\phi(p))( \la+|p|^{
2}/2)}\,dp\nonumber\\ && -\la\int_{[-\pi,\pi]^2}
{\phi(p)-1\over (1-e^{ -\la }\phi(p))( \la+|p|^{
2}/2)}\,dp\nonumber\\ && +(e^{ -\la }-1+\la)\int_{[-\pi,\pi]^2}
{\phi(p)\over (1-e^{ -\la }\phi(p))( \la+|p|^{ 2}/2)}\,dp\nonumber
\end{eqnarray}

Since
\begin{equation} |e^{ ip\cdot x}-1-ip\cdot x+(p\cdot x)^{
2}/2|\leq c(p\cdot x)^{  2+\de}\label{45.5}
\end{equation} for some $c<\ff$ we have by our assumptions that
\begin{equation} |\phi(p)-1+|p|^{ 2}/2|\leq c'|p|^{
2+\de}.\label{45.6}
\end{equation} this implies that
\begin{equation} |\phi(p)-1|\leq c''|p|^{ 2}\label{45.7}
\end{equation} for $p\in [-\pi,\pi]^2$ and \begin{equation}
1-e^{ -\la }\phi(p)\geq \bar{c}(\la+|p|^{ 2})\label{45.8}
\end{equation} for some $\bar{c}>0$ and sufficiently small $\la$.
Hence
\begin{eqnarray} && (e^{ -\la }-1+\la)\int_{[-\pi,\pi]^2}
{|\phi(p)|\over (1-e^{ -\la }\phi(p))( \la+|p|^{
2}/2)}\,dp\label{45.9}\\ &&
\leq c\la^{ 2} \int_{[-\pi,\pi]^2} {1\over ( \la+|p|^{ 2})^{
2}}\,dp\nonumber\\ &&
\leq c\la\int_{[-\pi/\sqrt{\la},\pi\sqrt{\la}]^2} {1\over ( 1+|p|^{
2})^{ 2}}\,dp=O( \la)\nonumber
\end{eqnarray} and
\begin{eqnarray} &&
\la\int_{[-\pi,\pi]^2} {|\phi(p)-1|\over (1-e^{ -\la }\phi(p))(
\la+|p|^{ 2}/2)}\,dp\label{45.10}\\ &&
\leq c\la\int_{[-\pi,\pi]^2} {|p|^{ 2}\over ( \la+|p|^{ 2})^{
2}}\,dp\nonumber\\ &&
\leq c\la\int_{[-\pi/\sqrt{\la},\pi\sqrt{\la}]^2} {|p|^{ 2}\over ( 
1+|p|^{ 2})^{ 2}}\,dp=O( \la\log ( 1/\la)).\nonumber
\end{eqnarray}

Setting $f( p)=\phi(p)-1+|p|^{ 2}/2$, and using (\ref{45.6}) we
see  that
\[\int_{[-\pi,\pi]^2}
     {|f( p)|\over |(1-\phi(p))|\,|p|^{ 2}/2}\,dp<\ff\]  Consider then
\begin{eqnarray} &&\qquad
\int_{[-\pi,\pi]^2} {f( p)\over (1-e^{ -\la }\phi(p))( \la+|p|^{
2}/2)}\,dp-\int_{[-\pi,\pi]^2}
     {f( p)\over (1-\phi(p))|p|^{ 2}/2}\,dp\label{45.20}\\
&&=\int_{[-\pi,\pi]^2} {f( p)\over (1-e^{ -\la }\phi(p))( \la+|p|^{
2}/2)}\,dp-\int_{[-\pi,\pi]^2}
     {f( p)\over (1-e^{ -\la }\phi(p))|p|^{ 2}/2}\,dp\nonumber\\
     &&+\int_{[-\pi,\pi]^2} {f( p)\over (1-e^{ -\la }\phi(p))|p|^{
2}/2}\,dp-\int_{[-\pi,\pi]^2}
     {f( p)\over (1-\phi(p))|p|^{ 2}/2}\,dp\nonumber
\end{eqnarray} We have
\begin{eqnarray} &&\qquad
\int_{[-\pi,\pi]^2} {f( p)\over (1-e^{ -\la }\phi(p))( \la+|p|^{
2}/2)}\,dp-\int_{[-\pi,\pi]^2}
     {f( p)\over (1-e^{ -\la }\phi(p))|p|^{ 2}/2}\,dp\label{45.21}\\
&& =-\int_{[-\pi,\pi]^2} {f( p)\la\over (1-e^{ -\la }\phi(p))(
\la+|p|^{ 2}/2)|p|^{ 2}/2}\,dp=O( \la^{ \de  }\log (
1/\la))\nonumber
\end{eqnarray} and the last line in (\ref{45.20}) can be bounded
similarly.  This completes the proof of Lemma \ref{lem-green}.
\qed

\begin{lemma}\label{lem-momgreen} Let $S_{n}$ be as above. For
all
$m\geq 1$
\be
     \|G_{ \la}\|_{ m}=O( \la^{-1/m }) \mbox{ as } \la\rar
0\label{45.31}
\ee and
\be
     \|G_{ \la}-G_{ \la'}\|_{ m}=O(|\la-\la'| (\sqrt{\la\la'})^{-1/m })
\mbox{ as } \la\rar 0.\label{45.31r}
\ee For all $m\geq 2$ and $z\in Z^{ 2}$
\be
     \|\Delta_{ z/\sqrt{\la}}G_{ \la}\|_{ m}\leq  c'
|z|^{2/m}\la^{-1/m }(\log (1/\la))^{ 1-1/m}\label{45.41a}
\ee
     and for any $0<\beta<1$
\be
     \|\Delta_{ z/\sqrt{\la}}G_{ \la}\|_{ m}\leq  c'
|z|^{\beta/m}\la^{-1/m }\label{45.41}
\ee and
\be
     \|\(\prod_{ i=1}^{ k}\Delta_{ z_{ i}/\sqrt{\la}}\)G_{ \la}\|_{
m}\leq c'
\(\prod_{ i=1}^{ k} |z_{ i}|^{\beta/mk}\)\la^{-1/m
}.\label{45.41m}
\ee
\end{lemma}

{\bf Proof of Lemma \ref{lem-momgreen}:}  By \cite{S}, p.~77, we
know that
$q_n(x)\leq c_1/n$, where $q_n$ is the transition probability for
$S_n$. So
\[ \| q_n\|_m^m =\sum_{x\in Z^2} q_n(x)^m \leq c_1^{m-1}
n^{-m+1}
\sum_{z\in Z^2} q_n(x)=c_1^{m-1} n^{-m+1}.
\] Then
\[
\| G_\la\|_m \leq \sum_{n=0}^\infty e^{-\la n} \|q_n\|_m.
\] Substituting the above estimate for $\|q_n\|_m$ and breaking
the sum into the sum over $n\leq 1/\la$ and the sum over
$n>1/\la$, we easily obtain (\ref{45.31}).

(\ref{45.31r}) follows from (\ref{45.31}) and the resolvent
equation
\be G_{ \la}-G_{ \la'}=(\la'-\la )G_{ \la}\ast\,G_{
\la'}.\label{45.31r2}
\ee

By Proposition 2.1 of \cite{BK2}, for each $\beta\in (0,1]$ there
exists a constant
$c_\beta$ such that
\[ |q_n(x)-q_n(y)|\leq c_\beta n^{-1} (|x-y|/\sqrt n)^\beta.
\] So for any fixed $w\in Z^2$
\begin{eqnarray}
\| q_n(\cdot+w)-q_n(\cdot)\|_m^m&\leq &
\|q_n(\cdot+w)-q_n\|_\infty^{m-1}
\sum_{x\in Z^2} (q_n(x+w)q_n(x)\nn\\ & \leq & 2(c_\beta
n^{-1}(|w|/\sqrt n)^\beta)^{m-1}.\nn
\end{eqnarray} We take $m$th roots, substitute into
\[\| G_\la(\cdot+w)-G_\la(\cdot)\|_m\leq
\sum_{n=0}^\infty e^{-\la n} \| q_n(\cdot+w)-q_n(\cdot)\|_m,
\] break the sum into the sum over $n\leq 1/\la$ and $n>1/\la$,
and let
$w=z/\sqrt \la$ to obtain (\ref{45.41}).

For (\ref{45.41m}) we note that for each $j$ we can write
$\(\prod_{ i=1}^{ k}\Delta_{ z_{ i}/\sqrt{\la}}\)G_{ \la}$ as a
sum of $2^{ k-1}$ terms of the form
$\Delta_{ z_{ j}/\sqrt{\la}}G_{ \la}( z+b)$ for some $b$ so that by
(\ref{45.41})
\be
     \|\(\prod_{ i=1}^{ k}\Delta_{ z_{ i}/\sqrt{\la}}\)G_{ \la}\|_{ 
m}\leq c' 2^{ k-1}|z_{ j}|^{\beta/m}\la^{-1/m }.\label{45.41p}
\ee
  We have inequality (\ref{45.41p}) for
$j=1, \ldots, k$. If we take the product of these $k$ inequalities
and then take $k$th roots, we have (\ref{45.41m}).
\qed

\def\noopsort#1{} \def\printfirst#1#2{#1}
\def\singleletter#1{#1}
      \def\switchargs#1#2{#2#1}
\def\bibsameauth{\leavevmode\vrule height .1ex
      depth 0pt width 2.3em\relax\,}
\makeatletter
\renewcommand{\@biblabel}[1]{\hfill#1.}\makeatother
\newcommand{\bysame}{\leavevmode\hbox to3em{\hrulefill}\,}

\end{document}